\documentclass{psp}

\ifprodtf \else
  \checkfont{msam10}
  \iffontfound
    \IfFileExists{amsfonts.sty}
      {\typeout{^^JFound AMS Symbol fonts on the system, using the
                'amsfonts' package.^^J}%
       \usepackage{amsfonts}%
      }
      {\providecommand\mathbb[1]{\mathsf{##1}}
       \providecommand\mathfrak[1]{\mathcal{##1}}}
  \else
    \providecommand\mathbb[1]{\mathsf{#1}}
    \providecommand\mathfrak[1]{\mathcal{#1}}
  \fi

\fi

\ifprodtf \else
  \IfFileExists{amsbsy.sty}
    {\typeout{^^JFound the 'amsbsy' package on the system, using it.^^J}%
     \usepackage{amsbsy}}
    {}

\fi

\newromanexpr\Hess{Hess}

\def\bC{{\mathbb{\overline{C}}}}
\def\ord{{\mathrm{ord}}}

\newcommand{\binom}[2]{\left(\begin{array}{c}#1\\#2\end{array}\right)}

\begin{document}

\title[Higher order Briot--Bouquet
differential equations]
{Meromorphic solutions of higher order Briot--Bouquet
differential equations}

\author[A. Eremenko, L.W. Liao and T.W. Ng]{ALEXANDRE E EREMENKO
\thanks{Supported by NSF grants
DMS-0555279 and DMS-0244547.\mbox{\hspace{4truecm}}}\\
Department of Mathematics, Purdue University,\\ West Lafayette IN, \textup{47907} USA\\
e-mail: {\rm eremenko@math.purdue.edu} \\
\\
{\rm LIANGWEN LIAO}\thanks{Partially supported
by the grant of the education department of Jiangsu Province, China \textup{07}KJB\textup{110069}.} \\
Department of Mathematics, Purdue University,\\Nanjing University,Nanjing, \textup{210093}, China\\
e-mail: {\rm maliao@nju.edu.cn}\\
\\
{\rm TUEN WAI NG}\thanks{Partially supported by RGC grant
HKU 7020/03P, and NSF grant DMS-0244547.\mbox{\hspace{0.9truecm}}}\\
Department of Mathematics, The University of Hong Kong,\\
Pokfulam, Hong Kong\\
e-mail: {\rm ntw@maths.hku.hk}}

\volume{121}
\pubyear{2009}
\setcounter{page}{1}
\maketitle

\begin{abstract}
For differential equations
$P(y^{(k)},y)=0,$ where $P$ is a polynomial,
we prove that all meromorphic solutions
having at least one pole
are elliptic functions, possibly degenerate.
\end{abstract}

\section{Introduction}

According to a theorem of Weierstrass, meromorphic functions
$y$ in the complex plane $\mathbb C$ that satisfy an algebraic addition theorem
\begin{equation}\label{1}
Q(y(z+\zeta),y(z),y(\zeta))\equiv 0,
\quad\mbox{where $Q\neq 0$ is a polynomial},
\end{equation}
are elliptic functions, possibly degenerate \cite{Ph,Akh}.

\newpage
More precisely,
let us denote by $W$ the class of meromorphic functions in $\mathbb C$
that consists of doubly periodic functions, rational functions
and functions of the form $R(e^{az})$ where $R$ is rational and $a\in\mathbb C$.
Then each function $y\in W$ satisfies an identity of the form
(\ref{1}), and conversely, every meromorphic
function\footnote{A ``meromorphic function''
in this paper means a function meromorphic in the complex plane,
unless some other domain is specified. See \cite{Ph,Ritt} for discussion
of the equation (\ref{1}) in more general classes of functions.}
that satisfies such
an identity belongs to $W$.

One way to prove this result is to differentiate (\ref{1}) with respect to
$\zeta$ and then set $\zeta=0$. Then we obtain a Briot--Bouquet
differential equation
$$P(y',y)=0.$$
The fact that every meromorphic solution of such an equation belongs to $W$
was known to Abel and Liouville, but probably it was stated for the
first time in the work of Briot and Bouquet \cite{BB,BB1}.

Here we consider meromorphic solutions of
higher order Briot--Bouquet equations
\begin{equation}\label{BB}
P(y^{(k)},y)=0,\quad\mbox{where $P$ is a polynomial}.
\end{equation}
Picard \cite{Picard} proved that for $k=2$, all meromorphic solutions
belong to the class $W$. This work was one of the first applications
of the famous Picard's theorems on omitted values.

In the end of 1970-s Hille \cite{H1,H2,H3,H4} considered
meromorphic solutions of (\ref{BB}) for arbitrary $k$.
The result of Picard was already forgotten,
and Hille stated it as a conjecture. Then Bank and Kaufman
\cite{BK} gave another proof of Picard's theorem.

These investigations were continued in \cite{E}. To state
the main results from \cite{E} we assume without loss of generality
that the polynomial $P$ in (\ref{BB}) is irreducible.
Let $F$ denote the compact Riemann
surface defined by the equation
\begin{equation}
\label{alg}
P(p,q)=0.
\end{equation}
Then every meromorphic solution $y$ of (\ref{BB}) defines a holomorphic
map $f:\mathbb C\to F$. According to another theorem of Picard,
a Riemann surface which admits a non-constant holomorphic
map from $\mathbb C$ has to be of genus $0$ or $1$,
(\cite{Picard1887}, see also \cite{BN}). The following theorems
were proved in \cite{E}:

\vspace{.1in}

\textsc{Theorem A.}
{\em If $F$ is of genus $1$, then every meromorphic
solution of $(\ref{BB})$ is an elliptic function.}
\vspace{.1in}

\textsc{Theorem B.}
{\em If $k$ is odd, then every meromorphic solution
of $(\ref{BB})$ having at least one pole, belongs to the class $W$.}
\vspace{.1in}

The main result of the present paper is the extension of
Theorem~B to the case of even $k$.
\vspace{.1in}

\textsc{Theorem 1.}
{\em If $y$ is a meromorphic solution of
an equation $(\ref{BB})$ and $y$ has at least one pole, then $y\in W$.}
\vspace{.1in}

This can be restated in the following way.
{\em Let $y$ be a meromorphic function in the plane
which is not entire and does not belong to $W$.
Then $y$ and $y^{(k)}$ are algebraically independent.}

It is easy to see that for every function $y$ of class $W$
and every natural integer $k$ there exists an
equation of the form (\ref{BB}) which $y$ satisfies.

It is not true that all meromorphic solutions of higher order
Briot--Bouquet equations belong to $W$, a simple counterexample
is $y^{\prime\prime\prime}=y$. We don't know whether
non-linear irreducible counterexamples exist.

In the process of proving of Theorem~1 we will establish  an estimate
of the degrees of possible meromorphic solutions in terms of
the polynomials $P$. Here by degree of a function of class $W$
we mean the degree of a rational function $y$, or the degree of $R$
in $y(z)=R(e^{az})$, or the number of poles in the fundamental parallelogram
of an elliptic function $y$.
Thus our result permits in principle the determination of all meromorphic
solutions having at least one pole of a given equation (\ref{BB}).

Our method of proof is based on the so-called ``finiteness property''
of certain autonomous differential equations:
there are only finitely many formal Laurent series
with a pole at zero that satisfy these equations.
The idea seems to occur for the first time in
\cite[p. 274]{H1} but the argument given there contains a mistake.
This mistake was corrected in \cite{E}.
Later the same method was applied in \cite{Hal} and \cite{E3}
to study meromorphic solutions of other differential equations.

\section{Preliminaries}

We will use the following refined version of Wiman--Valiron
theory which is due to Bergweiler, Rippon and Stallard.

Let $y$ be a meromorphic function and $G$ a component of
the set $\{ z:|y(z)|>M\}$ which contains no poles
(so $G$ is unbounded). Set
$$M(r)=M(r,G,y)=\max\{ |y(z)|:|z|=r,\, z\in G\},$$
and
\begin{equation}\label{a(r)}
a(r)=d\log M(r)/d\log r=rM'(r)/M(r).
\end{equation}
This derivative exists for all $r$ except possibly a discrete
set. According to a theorem of Fuchs \cite{F},
$$a(r)\to\infty,\quad r\to\infty,$$
unless the singularity of $y$ at $\infty$ is a pole.
For every $r>r_0=\inf\{|z|:z\in G\}$ we choose a point
$z_r$ with the properties $|z|=r,\; |y(z_r)|=M(r)$.
\vspace{.1in}

\textsc{Theorem C.} 
{\em For every $\tau>1/2$,
there exists a set $E\subset[r_0,+\infty)$ of
finite logarithmic measure, such that for $r\in[r_0,\infty)\backslash E$,
the disk
$$D_r=\{ z:|z-z_r|<ra^{-\tau}(r)\}$$
is contained in $G$ and
we have
\begin{equation}
\label{one}
y^{(k)}(z)=\left(\frac{a(r)}{z}\right)^k\left(\frac{z}{z_r}\right)^{a(r)}
y(z)(1+o(1)),\quad r\to\infty,\quad z\in D_r.
\end{equation}
}
\vspace{.1in}

When $y$ is entire, this is a classical theorem of Wiman.
Wiman's proof used power series, so it cannot be extended
to the situation when $y$ is not entire.
A more flexible proof, not using power series is due to Macintyre \cite{Mc};
it applies, for example to functions analytic and unbounded in $|z|>r_0$.
The final result stated above was recently established in
\cite{Be}.

\section{Proof of Theorem 1}
\vspace{.1in}

In what follows, we always assume that the polynomial $P$ in (\ref{BB})
is irreducible.

To state a result of \cite{E} which we will need,
we introduce the following
notation. Let $A$ be the field of meromorphic functions on $F$.
The elements of $A$ can be represented as rational functions
$R(p,q)$ whose denominators are co-prime with $P$.
In particular,  $p$ and $q$ in (\ref{alg}) are elements of $A$.
For $\alpha\in A$ and a point $x\in F$, we denote by $\ord_x\alpha$
the order of $\alpha$ at the point $x$.
Thus if $\alpha(x)=0$ then $\ord_x\alpha$ is the multiplicity
of the zero $x$ of $\alpha$, if $\alpha(x)=\infty$ then
$-\ord_x\alpha$ is the multiplicity of the pole, and
$\ord_x\alpha=0$ at all other points $x\in F$.

Let $I\subset F$ be the set of poles of $q$.
For $x\in I$ we set $\kappa(x)=\ord_xp/\ord_xq.$

\vspace{.1in}

\textsc{Theorem D.} 
{\em Suppose that an irreducible equation
$(\ref{BB})$ has a transcendental meromorphic solution $y$.
Let $f:\mathbb C\to F$ be the holomorphic map defined by
$z\mapsto(y^{(k)}(z),y(z))$.
Then:
\newline
a) Every pole of $p$ belongs to $I$.
\newline
b) For every $x\in I$, the number $\kappa(x)$ is either $1$ or
$1+k/n$, where $n$ is a positive integer.
\newline
c) If $\kappa(x)=1+k/n$ for some $x\in I$, then the equation
$f(z)=x$ has infinitely many solutions, and all these solutions
are poles of order $n$ of $y$.
\newline
d) If $\kappa(x)=1$ for some $x\in I$, then the equation
$f(z)=x$ has no solutions.}
\vspace{.1in}

Picard's theorem on omitted values implies that $\kappa(x)=1$ can
happen for at most two points $x\in I$. For the convenience
of
the reader we include a proof of Theorem~D in the Appendix.

The numbers $\kappa(x)$ can be easily determined from the Newton
polygon of $P$. Thus Theorem D gives several effective
necessary conditions
for the equation (\ref{BB}) to have meromorphic or entire solutions.
\vspace{.1in}

\noindent
{\em Remark.} The proof of Theorem~D in \cite{E} uses Theorem~C
which was stated in \cite{E} but not proved.
One can also give an
alternative proof of Theorem~D, using Nevanlinna theory
instead of Theorem~C, by the arguments similar to
those in \cite{E2}.

\vspace{.1in}

\textsc{Lemma 1.} {\em Suppose that $y$ is a meromorphic solution
of $(\ref{BB})$. If $\kappa(x)=1$ for some $x\in I$ then
$y$ has order one, normal type.}
\vspace{.1in}

{\em Proof.} In view of Theorem A and Theorem D, d), we conclude that
the genus of $F$ is zero. Therefore, we can find $t=R(p,q)$ in $A$
which has a single simple pole at $x$. Then $w=R(y^{(k)},y)$ is an entire
function by Theorem D, d). As $t$ has a simple pole at $x$, the element $1/t\in A$ is
a local parameter at $x$, and
 in a neighborhood of $x$ we have
$$q=at^m+\ldots\quad\mbox{and}\quad p=bt^m+\ldots,$$
where $-m=\ord_xp=\ord_xq$ as $\kappa(x)=1$, and the dots stand for the terms
of degree smaller than $m$. Substituting $p=y^{(k)}$ and
$q=y$ and differentiating the first equation $k$ times we obtain
for $w$ a differential equation of the form
\begin{equation}
\label{two}
\frac{d^k}{dz^k}w^m +\cdots=(b/a)w^m,
\end{equation}
where the dots stand for the terms of degree smaller than $m$.
Now we use a standard
argument of Wiman--Valiron theory. Applying Theorem~C
to the entire function $w^m$, with
$G=\mathbb C$ and $z=z_r$, we compare the asymptotic
relations (\ref{one}) and (\ref{two}) to conclude that
$a(r)\sim c r,$ where $c\neq 0$ is a constant.
This implies $\log M(r)\sim cr$, which means that 
$w$ is of order $1$, normal type.
So $y$ is also of order $1$, normal type,
because $w$ and $y$ satisfy a polynomial relation
of the form $P(y,w)=0$, where $P$ is a polynomial
with constant coefficients.
\vspace{.1in}

\textsc{Lemma 2.} {\em Suppose that $y$ is a meromorphic solution
of $(\ref{BB})$. If $\kappa(x_1)=\kappa(x_2)=1$ for two different
points $x_1$ and $x_2$ in $I$, then $y$ is 
a rational function of $e^{az}$, where $a\in\mathbb C$.}
\vspace{.1in}

{\em Proof}. As in the previous lemma, the genus
of $F$ is zero. Let $t=R(p,q)$ be a function in $A$ with
a single simple pole at $x_1$ and a single simple zero
at $x_2$. Then $w=R(y^{(k)},y)$ is an entire function
of order $1$, normal type (by Lemma~1) omitting $0$ and $\infty$
(by Theorem D, d). So $w(z)=e^{az}$ for some $a\in\mathbb C$. 
Since $t$ is a generator of $A$, by L\"uroth's theorem, 
both $p$ and $q$ are rational functions of $t$ 
and the lemma follows.

\vspace{.1in}

\textsc{Lemma 3.} {\em Suppose that $k$ is even,
the Riemann surface $F$ is
of genus zero, $y$ is a non-constant
 meromorphic solution of $(\ref{BB})$,
and $\kappa(x)=1$ for at most one point $x\in I$.
Then the Abelian differential $pdq$ is exact,
that is $pdq=ds$ for some $s\in A$.}

\vspace{.1in}

{\em Proof.} It is sufficient to show that under the assumptions
of Lemma 3, the integral of $pdq$ over every closed path in $F$
is zero. As $F$ is of genus zero, we only have to consider
residues of $pdq$.
By Theorem D, a), all poles of our differential
belong to the set $I$.

Consider first a point $x\in I$ with $\kappa(x)=1+k/n$.
By Theorem D, c),  we have a meromorphic solution $y$ with
a pole of order $n$ at zero, such
that the corresponding function $f$ has the property
$f(0)=x$. In a neighborhood of $x$ we have a Puiseaux
expansion
$$pdq=\sum_{j=J}^\infty c_jq^{-j/m}dq$$
with some positive integer $m$.
We substitute $p=y^{(k)},\; q=y$ and obtain
\begin{equation}
\label{j}
y^{(k)}y'=\sum_{j\neq -m} c_jy^{-j/m}y' +ry^{-1}y',
\end{equation}
where $r=c_m$ is the residue of $pdq$ at $x$.
Now we notice that for even $k$,
\begin{equation}
\label{integration}
y^{(k)}y'=\frac{d}{dz}\left\{ y^{(k-1)}y'-y^{(k-2)}y^{\prime\prime}+\ldots
\pm\frac{1}{2} (y^{(k/2)})^2\right\}.
\end{equation}
Using this, we integrate (\ref{j}) over a small circle around $0$
in the $z$-plane, described $m$ times anticlockwise. We obtain
that $2\pi i mr=0$, so $r=0$.

Now we consider a point $x\in I$ with $\kappa(x)=1$.
By the assumptions of the lemma, there is at most one such point.
Then the residue of $pdq$ at $x$ is zero because the sum of all
residues of a differential on a compact Riemann surface is zero.
This proves the lemma.
\vspace{.1in}

Using (\ref{integration}) and Lemma 3, if the assumptions of Lemma 3
are satisfied, we can rewrite our differential equation
\begin{equation}\label{ish}
y^{(k)}=p(y)
\end{equation}
 as
\begin{equation}
\label{rewrite}
y^{(k-1)}y'-y^{(k-2)}y^{\prime\prime}+\ldots
\pm \frac{1}{2}(y^{(k/2)})^2=s(y)+c,
\end{equation}
where $s\in A$ is an integral of the exact differential $pdq$,
and $c$ is a constant that depends on the particular solution $y$.
We have the relation $p(y)=ds/dy$.
\vspace{.1in}

\textsc{Lemma 4.} {\em For a given differential equation
of the form $(\ref{rewrite})$, there are only finitely many formal Laurent
series with a pole at zero that satisfy the equation.}
\vspace{.1in}

{\em Proof.} By making a linear change of the independent
variable, we may assume that
$$s(y)=y^{2+k/n}+\ldots.$$
Then $$p(y)=(2+k/n)y^{1+k/n}+\ldots .$$
Now we substitute a Laurent series with undetermined
coefficients
\begin{equation}
\label{19}
y(z)=\sum_{j=0}^\infty c_jz^{-n+j}
\end{equation}
to the equation (\ref{ish}), which is a consequence of (\ref{rewrite}).
With even $k$ we have:
\begin{eqnarray*}
y^{(k)}(z)&=&\frac{(k+n-1)!}{(n-1)!}c_0z^{-n-k}+
\frac{(k+n-2)!}{(n-2)!}c_1z^{-n-k-1}\\
&&+\ldots+k!c_{n-1}z^{-k-1}\\
&&+k!c_{n+k}+\frac{(k+1)!}{1!}c_{n+k+1}z+\frac{(k+1)!}{2!}c_{n+k+2}z^2
+\ldots;
\end{eqnarray*}
and
\begin{eqnarray*}
y^{1+k/n}(z)&=&z^{-k-n}\left[c_0^{1+k/n}+\left((1+k/n)c_0^{k/n}c_1+
(\ldots)_1\right)z\right.\\
&&+\left((1+k/n)c_0^{k/n}c_2+(\ldots)_2\right)z^2+\ldots\\
&&\left.+\left((1+k/n)c_0^{k/n}c_j+(\ldots)_j\right)z^j+\ldots\right].
\end{eqnarray*}
In the last formula, the symbol $(\ldots)_j$ stands for a finite sum of
products of the coefficients of the series (\ref{19}) which contain no
coefficients $c_i$ with $i 
\geq
 j$. Substituting to (\ref{ish}) and comparing
the coefficients at $z^{-k-n}$ we obtain
$$\frac{(k+n-1)!}{(n-1)!}c_0=(2+k/n)c_0^{1+k/n}.$$
This equation has finitely many non-zero roots $c_0$.
We have
\begin{equation}
\label{21}
(2+k/n)c_0^{k/n}=\frac{(k+n-1)!}{(n-1)!}.
\end{equation}
Further we obtain
\begin{equation}
\label{22}
\frac{(k+n-2)!}{(n-1)!}c_1=(2+k/n)c_0^{k/n}(1+k/n)c_1+(\ldots)_1.
\end{equation}
Substituting here the value of $(2+k/n)c_0^{k/n}$ from (\ref{21}),
we see that the coefficient at $c_1$ is different from zero, because
$$\frac{(k+n-2)!}{(n-2)!}\neq\frac{(k+n-1)!}{(n-1)!}
\frac{k+n}{n}.$$
Thus $c_1$ is uniquely determined from (\ref{22}). The situation
is analogous for all coefficients $c_j$ with $j<n+k$. These
coefficients are uniquely determined from the equation (\ref{ish})
once $c_0$ is chosen.

Now we consider the coefficients $c_{n+k+j}$ with $j\geq 0$.
We have
$$\frac{(k+j)!}{j!}c_{n+k+j}=(2+k/n)c_0^{k/n}\frac{n+k}{n}c_{n+k+j}+
(\ldots)_{n+k+j}.$$
Again we substitute the value of $(2+k/n)c_0^{k/n}$ from (\ref{21}) and
conclude that the coefficient at $c_{n+k+j}$ equals
$$\frac{(k+j)!}{j!}-\frac{(k+n)!}{n!}.$$
This coefficient is zero for a single value
of $j$, namely $j=n$. Thus $c_{2n+k}$ cannot be determined from
the equation (\ref{ish}), but once $c_0$ and $c_{2n+k}$ are chosen,
the rest
of the coefficients of the series (\ref{19}) are determined uniquely.

To determine $c_{2n+k}$ we invoke the equation (\ref{rewrite}):
\begin{equation}
\label{44}
y^{(k-1)}y'-y^{(k-2)}y^{\prime\prime}+\ldots
\pm \frac{1}{2}(y^{(k/2)})^2=y^{2+k/n}+\ldots,
\end{equation}
where the dots stand for the terms of lower degrees.
We have
\begin{eqnarray*}
y'(z)&=&-nc_0z^{-n-1}+\ldots+c_{2n+k}(n+k)z^{n+k-1}+\ldots,\\
y^{\prime\prime}&=&n(n+1)c_{0}z^{-n-2}+\ldots+c_{2n+k}(n+k)(n+k-1)z^{n+k-2}
+\ldots,\\
\ldots&&\ldots,\\
y^{(k-1)}&=&-n(n+1)\ldots(n+k-2)c_0z^{-n-k+1}+\ldots\\
&&+c_{2n+k}(n+k)(n+k-1)\ldots(n+2)z^{n+1}+\ldots .
\end{eqnarray*}
Substituting this to our equation (\ref{44}) we write the condition
that the constant terms in both sides of (\ref{44}) are equal.
This condition is a polynomial equation in $c,c_0,\ldots,c_{2n+k}$
(it is linear with respect to $c_{2n+k}$)
and the coefficient at $c_{2n+k}$ in this equation
equals
$$c_0\sum_{m=0}^{k-1}\frac{(n+m)!(n+k)!}{(n+m+1)!(n-1)!}.$$
This expression is not zero because
each term of the sum is positive.
Thus $c_{2n+k}$ is determined uniquely,
and this completes the proof of the lemma.
\vspace{.1in}

{\em Remark.} It follows from this proof that the only
meromorphic solutions of the differential equations
$$y^{(k)}=y^m$$
are exponential polynomials when $m=1$ and functions
$c(z-z_0)^{-n}$ where $m=1+k/n,\; z_0\in\mathbb C$ and $c$ is
an appropriate constant.
\vspace{.1in}

The rest of the proof of Theorem 1 is a repetition of the argument from
\cite{E}.

By Theorems A and B, we may assume that $F$ is of genus zero,
and $k$ is even.
In view of Lemmas 2 and 3, it is enough to consider the case
that the differential $pdq$ is exact. Then every solution
of (\ref{BB}) also satisfies (\ref{rewrite}) with some constant $c$.

Assume that $y$ is a transcendental meromorphic solution of
(\ref{rewrite}), having at least one pole.
By Theorem D, d), c), $y$ has infinitely many poles $z_j,\; j=1,2,3,\ldots.$
The functions $y(z-z_j)$ satisfy the assumptions of Lemma 4, therefore
some of them are equal. We conclude that $y$ is a periodic function.
By making a linear change of the independent variable we may assume
that the smallest period is $2\pi i$.

Consider the strip $D=\{ z:0\leq\Im z<2\pi\}$.
\vspace{.1in}

{\em Case 1.} $y$ has infinitely many poles in $D$. Applying Lemma~4
again, we conclude that $y$ has a period in $D$, so $y$ is doubly periodic.
\vspace{.1in}

{\em Case 2.} $y$ is bounded in $D\cap\{ z:|\Re z|>C\}$ for some $C>0$.
Since $y$ is $2\pi i$-periodic, we have $y(z)=R(e^z)$ where
$R$ is meromorphic in $\mathbb C^*$. As $R$ is bounded in some neighborhoods
of $0$ and $\infty$, we conclude that $R$ is rational.
\vspace{.1in}

{\em Case 3.} $y$ has finitely many poles in $D$ and is unbounded
in $D\cap\{ z:|\Re z|>C\}$ for every $C>0$.
As $y$ is $2\pi i$-periodic, we write $y=R(e^z)$
where $R$ is meromorphic in $\mathbb C^*$. Now $R$ has finitely many
poles and is unbounded either in a neighborhood of $0$
or in a neighborhood of $\infty$. Suppose that it is
unbounded in a neighborhood of $\infty$.
Then the set  $\{ z:|R(z)|>M\}$, where $M$ is large enough
has an unbounded component $G$
containing no poles of $R$. On this component $G$, the function $R$
satisfies a differential equation
$$\sum_{m=1}^k\binom{k}{m}w^m\frac{d^mR}{dw^m}=(c+o(1))R^\kappa,$$
where $c$ is some constant and $\kappa=1$ or $\kappa$
is one of the numbers $1+k/n$ from Theorem D.
Applying Theorem~C in $G$ 
as we did in the proof of Lemma~1, we obtain that
$\kappa=1$ and that $R$ has a pole at infinity.
Similar argument works for the singularity at $0$,
so $R$ is rational, and this completes the proof.
\vspace{.1in}

\section{Appendix}

{\em Proof of Theorem D.} Statement a) is a special case
of \cite[Th. 10]{E2}, but we give a simple
independent proof using Theorem~C.
Proving it
by contradiction, suppose that
$p$ has a pole at a point $x\in F$ such that
$q(x)=b\in\mathbb C$. Let $D_\epsilon\subset \mathbb C$ be a disk
of radius $\epsilon$
centered at $b$, and $V_\epsilon\subset F$ a component
of $q^{-1}(D_\epsilon)$ containing $x$.
We assume that the disk
$D_\epsilon$ is so small that $V_\epsilon$
contains no other poles of $p$, except the pole at $x$.
Let $y$ be a meromorphic solution of our equation
(\ref{BB}) and consider the map $f:\mathbb C\to F$ given by
$f(z)=(y^{(k)}(z),y(z))$. The image of this map is dense
in $F$ and the point $x$ is evidently omitted
by $f$. Let $G_\epsilon\subset \mathbb C$ be a component of
the preimage $f^{-1}(D_\epsilon).$
Consider the meromorphic function
$w=1/(y-a)$. It is holomorphic and unbounded in $G_\epsilon$,
and $|w(z)|=1/\epsilon$ for $z\in\partial G_\epsilon$.
We conclude that
$G_\epsilon$ is unbounded.
Now we apply Theorem C to $w$ in $G_\epsilon$.

Set $M(r)=\max\{|w(z)|:|z|=r, z\in G_\epsilon\}$ 
and let $a(r)$ be defined as in (\ref{a(r)}).
For any $r>r_0=\inf\{|z|:z\in G_\epsilon\}$, we choose
a point $z_r$ with $|z|=r$ and $|w(z_r)|=M(r)$. 
By Theorem~C, we have 
\begin{equation}\label{e3}
|w^{(j)}(z_r)|=\left(\frac{a(r)}{r}\right)^j
|w(z_r)|(1+o(1))=\frac{a(r)^j}{r^j} M(r)(1+o(1))
\end{equation}
where $r\rightarrow\infty$ outside a
set of finite logarithmic measure.

{}From Lemma 6.10 of \cite{Be}, we have for every $\beta>0$,
\begin{equation}\label{e4}
(a(r))^\beta=o(M(r)),  
\end{equation}
as $r\rightarrow\infty$ outside a set
of finite logarithmic measure. 

Differentiating the equation $y=1/w+a$ we obtain
\begin{equation}
\label{raz}
y^{(k)}=\frac{1}{w} Q\left(\frac{w'}{w},
\frac{w''}{w},\cdots,\frac{w^{(k)}}{w}\right),\end{equation}
where $Q$ is a polynomial.
On the other hand, from the Puiseaux expansion at the
point $x$ we obtain
\begin{equation}
\label{dwa}
y^{(k)}=(c+o(1))w^\alpha,\quad w\to\infty,
\end{equation}
where $c\neq 0$ is a constant and $\alpha>0$.
Combining (\ref{raz}) and (\ref{dwa})
we obtain
$$Q\left(\frac{w'}{w},\ldots,\frac{w^{(k)}}{w}\right)
=(c+o(1))w^{1+\alpha}.$$
Inserting to this asymptotic relation $z=z_r$ and
using (\ref{e3}) and (\ref{e4}) we obtain a contradiction
which proves a).
\vspace{.1in}

Consider now a point $x\in I$. From the Puiseaux expansion
we obtain
\begin{equation}
\label{pui}
y^{(k)}=(c+o(1))y^{\kappa(x)},\quad y\to\infty.
\end{equation}
If $x$ has a preimage under the map $f$, then this preimage
is a pole $z_0$ of $y$. If this pole is of order $n$
we have $y(z)\sim c_1(z-z_0)^{-n}$ and $y^{(k)}(z)\sim
c_2(z-z_0)^{-n-k}$ as $z\to z_0$. Substituting
to (\ref{pui}) we conclude that $\kappa(x)=1+k/n$.
Thus if $x$ has at least one preimage under $f$ then
$\kappa(x)=1+k/n$ with a positive integer $n$, and every
preimage of $x$ is a pole of order $n$ of $y$.
This implies d).

Now suppose that a point $x\in I$ has only finitely
many preimages. Let $U_\epsilon=\{ z\in \bC:|z|>1/\epsilon
\}$
be a neighborhood
of infinity, and $V_\epsilon\subset F$ a component
of the preimage $q^{-1}(U_\epsilon)$.
We may assume that $\epsilon>0$ is so small
that $V_\epsilon$ does not contain other poles of $q$
except $x$. Let $G_\epsilon$
be a component of the preimage
$f^{-1}(V_\epsilon)$. If $G_\epsilon$ is bounded
then $f:G_\epsilon\to U_\epsilon$
is a ramified covering of
a finite degree, and $f$ takes the value $x$ somewhere
in $G$.
As we assume that $f$ is transcendental but $x$ has
only finitely many preimages, there should exist
an unbounded component $G_\epsilon$. Choosing a smaller
$\epsilon$ if necessary, we achieve that this
unbounded component $G_\epsilon$
contains no $f$-preimages
of $x$.
Then $y$ is a holomorphic function in $G_\epsilon$,
$|y(z)|=1/\epsilon,\; z\in\partial G_\epsilon$,
and $y$ is unbounded
in $G_\epsilon$. Applying Theorem~C to the
function $y$ in $G_\epsilon$ we obtain
the asymptotic relation (\ref{one}). Putting $z=z_r$
in this relation, taking (\ref{e4}) into account,
and comparing with (\ref{pui}) we conclude that
$\kappa=1$ in (\ref{pui}). This implies c).
Thus in any case $\kappa=1+k/n$ or
$\kappa=1$, which proves
b).

\enddocument